\documentclass[reqno]{amsart}
\usepackage[all,cmtip]{xy}
\usepackage{hyperref}
\usepackage[latin1]{inputenc}
\usepackage{accents}

\newcommand\R{\mathbb R}

\newcommand\m{\mathbf m}
\newcommand\LL{\mathcal L}

\newcommand\OO{\mathcal O}
\renewcommand\S{\mathcal S}

\newcommand\PP{\mathbb P}
\renewcommand\AA{\mathbb A}

\newcommand\X{\mathcal X}
\newcommand\Y{\mathcal Y}
\newcommand\Z{\mathcal Z}
\newcommand\HH{\mathcal H}

\newcommand{\mult}{\operatorname{mult}}

\newcommand{\Spec}{\operatorname{Spec}\,}
\newcommand{\Hilb}{\operatorname{Hilb}}
\newcommand{\Sing}{\mathcal{I}}
\newcommand{\USing}{\mathcal{P}} 

\newcommand{\oover}[1]{\accentset{\circ}{#1}}

\makeatletter \@addtoreset{equation}{section} \makeatother

\newtheorem{thm}[equation]{Theorem}
\newtheorem{lem}[equation]{Lemma}
\newtheorem{cor}[equation]{Corollary}

\newtheorem{conj}[equation]{Conjecture}

\theoremstyle{definition}

\newtheorem{rmk}[equation]{Remark}

\title{An inequality between multipoint Seshadri constants}
\author{J. Roé}
\thanks{J. Roé was partially supported by MEC (Spain) grants
  MTM2005­01518 and MTM2006­11391, and by the CIRIT (Catalonia) grant
  2005SGR00787}
\author{J. Ross}
\date{\today}

\begin{document}
\bibliographystyle{abbrv}

\begin{abstract}
  Let $X$ be a projective variety of dimension $n$ and $L$ be a nef divisor on
  $X$.  Denote by $\epsilon_d(r;X,L)$ the $d$-dimensional Seshadri constant of
  $r$ very general points in $X$.  We prove that
  $$\epsilon_d(rs;X,L)\ge \epsilon_d(r;X,L)\cdot
  \epsilon_d(s;\PP^n,\OO_{\PP^n}(1)) \quad \text{ for } r,s\ge 1.$$ 
\end{abstract}

\maketitle

\section{Introduction}

Let $L$ be a nef divisor on a projective variety $X$ of dimension $n$. 
The Seshadri constant of dimension $d\le n$ at $r$ points $p_1,\ldots,p_r$ in 
the smooth locus of $X$
is defined (see \cite{szemberg(01):global_local_pos_lin_bdl}) to be the real number
\begin{equation*}
\epsilon_d(p_1,\ldots,p_r;X,L) =\inf\left\{ \left( \dfrac{L^d\cdot Z}{\sum \mult_{p_i}Z} \right) ^{1/d}: 
\begin{matrix}
Z \subset X \text{ effective cycle} \\
\text{of dimension }d
\end{matrix}
\right\}.
\end{equation*}

By semicontinuity of multiplicities if the points are in very general position then
the Seshadri constant does not depend on the actual points chosen
\eqref{lem:semicontinuous}.  Thus one can define
$$\epsilon_d(r;X,L) = \epsilon_d(p_1,\ldots,p_r;X,L)$$
where $p_1,\ldots,p_r$ is any collection of $r$ very general points in
$X$.  We shall prove the following inequality comparing Seshadri
constants of very general points $X$ and those in projective space.

\begin{thm}\label{thm:main}
  Let $X$ be an $n$ dimensional projective variety, $L$ be a nef
  divisor on $X$ and $r,s,d\ge 1$ be integers, with $d\le n$.  Then
$$\epsilon_d(rs;X,L) \ge \epsilon_d(r;X,L) \cdot \epsilon_d(s;\PP^n,\OO_{\PP^n}(1)).$$
\end{thm}

\begin{rmk}  When $d=1$ the Seshadri constant $\epsilon_1$ reduces to the usual Seshadri constant as introduced by Demailly \cite{demailly(92):sing_herm_metrics}, and in this case Theorem \ref{thm:main} is due to Biran \cite{biran(99):const_new_ample_divis_out_old_ones}.   When $r=1$, $d=n-1$ the Theorem is due to Ro\'e \cite{roe(04):relat_between_one_point_multi}.
\end{rmk}

It is well known and not hard to see that
$\epsilon_d(r;X,L)\le\sqrt[n]{L^n/r}$
for $d=1$ and $d=n-1$ (see remark \ref{rmk:grow}). On the other hand,
explicit values or even lower bounds are in general
hard to compute. Thus Theorem \ref{thm:main} should be seen as a lower
bound on $\epsilon_d(rs;X,L)$, which will be useful if Seshadri
constants on projective space are known. As an example, let us recall
the most general form of a famous conjecture by Nagata:

\begin{conj}[Nagata-Biran-Szemberg, \cite{szemberg(01):global_local_pos_lin_bdl}]
Let $X$ be an $n$ dimensional projective variety, $L$ be a nef
  divisor on $X$ and $1 \le d\le n$ be an integer.  Then there is a
  positive integer $r_0$ such that for every $r\ge r_0$,
$\epsilon_d(r;X,L) = \sqrt[n]{L^n/r}.$
\end{conj}

\begin{rmk}\label{rmk:grow}
  Already Demailly observed that $\epsilon_1(r;X,L)\le\epsilon_d(r;X,L)$ for all
  $d=2, \dots, n$ (which combined with the obvious equality
  $\epsilon_n(r;X,L)= \sqrt[n]{L^n/r}$ gives the upper bound
  $\epsilon_1(r;X,L)\le\sqrt[n]{L^n/r}$). Standard arguments also show
  that $\epsilon_{n-1}(r;X,L)\le\sqrt[n]{L^n/r}=\epsilon_n(r;X,L)$ and
  that $\epsilon_1(r;X,L) = \sqrt[n]{L^n/r}$ implies
  $\epsilon_d(r;X,L) = \sqrt[n]{L^n/r}$ for all $d=2, \dots n$, so if
  the conjecture above holds for $d=1$ then it holds for all $d$. In
  view of these facts, it is tempting to ask if the inequalities
  $\epsilon_{d_1}(r;X,L)\le\epsilon_{d_2}(r;X,L)$ hold for all $d_1\le d_2$.
\end{rmk}

In this context Theorem \ref{thm:main} implies, for instance,
that if (1) the Nagata-Biran-Szemberg conjecture is true for
$d$-dimensional Seshadri constants on $\PP^n$,
and (2) $\epsilon_d(1;X,L)=\sqrt[n]{L^n}$, then the
Nagata-Biran-Szemberg conjecture is true for the $d$-dimensional Seshadri
constants of $(X,L)$. See \cite{roe(04):relat_between_one_point_multi}
for more applications along this line.

\smallskip

The proof uses the degeneration to the normal cone of $r$ very general
points of $X$.  That is, let $p_1,\ldots,p_r$ be very general
points in $X$ and $\pi\colon \X\to \AA^1$ be the blowup of $X\times
\AA^1$ at points $p_i\times \{0\}$ for $1\le i\le r$.  The central
fibre of $\X$ over $0\in \AA^1$ is reducible, having one component
which is the blowup of $X$ at $p_1,\ldots,p_r$ (with exceptional
divisor $E_i\simeq \PP^{n-1}$ over $p_i$) and $r$ exceptional
components $F_i\simeq \PP^n$, with each $E_i$ glued to $F_i$ along a
hyperplane.

Regard a collection of $rs$ points in $X$ as coming from $r$
sub-collections each consisting of $s$ points.  Shrinking $\AA^1$ if
necessary for each $1\le i\le r$ we can find $s$ sections of $\pi$ that
pass through very general points in $F_i\backslash E_i$.  Then the inequality in
Theorem \ref{thm:main} follows by showing that the existence of highly 
singular cycles at $rs$ points of $X$ implies the existence of highly 
singular cycles at $rs$ points of the central fiber chosen at will, 
which means cycles with higher multiplicities at the $s$ chosen points.

It is clear that this argument actually produces something stronger.
For instance we do not have to divide the collection of $rs$ points
evenly:

\begin{thm}\label{thm:main2}
  Let $X$ be an $n$ dimensional projective variety and $L$ be a nef divisor on
  $X$.  Suppose $1\le s_1\le s_2\le \cdots \le s_r$ are integers and
  set $s= \sum_{i=1}^r s_i$.  Then
$$\epsilon_d\left(s;X,L\right) \ge \epsilon_d(r;X,L) \cdot \epsilon_d(s_r;\PP^n,\OO_{\PP^n}(1)).$$
\end{thm}

We can also consider weighted Seshadri constants as in
\cite{harbourne-roe(08):disc_beh_sesh_consts}, defined as follows.
In addition to the nef divisor $L$ on $X$ and $r$ points $p_1,\ldots,p_r$ in 
the smooth locus of $X$, let a nonzero real vector
${\bf \ell}=(l_1,\cdots,l_r)\in \mathbb R_{+}^{r}$ be given with each $l_i\ge 0$.   The Seshadri constant of dimension $d\le n$ at these points with weights ${\bf \ell}$
is defined to be the real number
$$\epsilon_d(l_1p_1, \dots, l_rp_r;X,L)
=\inf\left\{ \left( \dfrac{L^d\cdot Z}{\sum l_i\mult_{p_i}Z} \right) ^{1/d}: 
\begin{matrix}
Z \subset X \text{ effective cycle} \\
\text{of dimension }d
\end{matrix}
\right\}.
$$

Again by semicontinuity if the points are in very general position then
the weighted Seshadri constant does not depend on the actual points chosen
\eqref{lem:semicontinuous}.  Thus one can define
$$\epsilon(r,{\bf \ell};X,L) = \epsilon(l_1p_1,\ldots,l_rp_r;X,L)$$
where $p_1,\ldots,p_r$ is any collection of $r$ very general points in
$X$.

\begin{thm}\label{thm:glueing}
  Let $X$ be an $n$ dimensional projective variety, let $L$ be a nef
  divisor on $X$, and let $r, d \ge 1$  
be integers with $d \le n$. Suppose $1\le s_1\le s_2\le \cdots \le s_r$ are integers and
set $s= \sum_{i=1}^r s_i$. For each $i$ let ${\bf \ell}_i$ be a vector in $\mathbb R_{+}^{s_i}$ and set ${\bf \ell} = ({\bf \ell}_1, {\bf \ell}_2, \ldots,{\bf \ell}_r)\in \mathbb R_{+}^{s}$.
Then
$$\epsilon_d(s,{\bf \ell};X,L) \ge \epsilon_d(r;X,L) \cdot 
\inf_{i=1,\ldots,r}\epsilon_d(s_i,{\bf \ell}_i;\PP^n,\OO_{\PP^n}(1)).$$
\end{thm}

\begin{rmk}
Since the $L$-degree and multiplicity at a point of a $d$-dimensional
scheme coincide with those of the associated $d$-dimensional cycle, the
Seshadri constant can be equivalently defined as
\begin{equation*}
\epsilon_d(p_1,\ldots,p_r;X,L) =\inf\left\{ \left( \dfrac{L^d\cdot Z}{\sum \mult_{p_i}Z} \right) ^{1/d}: 
\begin{matrix}
Z \subset X \text{ closed subscheme} \\
\text{of dimension }d
\end{matrix}
\right\}.
\end{equation*}
For convenience, we shall work with this definition.
\end{rmk}

\vspace{4mm}{\bf Conventions: } 
By a variety
$X$ we mean a possibly reducible reduced scheme of finite type
over an algebraically closed uncountable field $K$.

If $X$ is a variety then a very general point $p\in X$ is a
point that lies outside a countable collection of proper subvarieties
of $X$.  We say that a collection of $r$ points $p_1,\ldots,p_r$ is
very general if the point $(p_1,\ldots,p_r)\in X^r$ in the $r$-th
power of $X$ is very general. The uncountable hypotheses on the base
field guarantees that a set of very general points is nonempty and
dense; over a countable or finite field claims on very general points
are void and $\epsilon(r;X,L)$ is not well defined.

\section{Preliminaries}

\subsection{Semicontinuity of Seshadri constants}

Our proof of Theorem \ref{thm:main} will rely on the 
semicontinuity property according to which Seshadri constants
can only ``jump down'' in families.   When $d=1$ this semicontinuity
is well-known and comes from interpreting Seshadri constants in terms of
ampleness of certain line bundles on a blowup of $X$
\cite[Thm. 5.1.1]{lazarsfeld(04):posit_in_algeb_geomet}.    To deal
with the general case $d\ge 1$, let $p\colon \X\to B$ be a projective
flat morphism, and for $b\in B$, denote $\X_b=p^{-1}(b)$. 
We assume that $\X$ is a variety, $B$ is a reduced and
irreducible scheme, and all the fibres $\X_b$ are
(possibly reducible)
varieties. The following result, well known in singularity theory,
follows from the Hilbert-Samuel stratification 
\cite{lejeune-jalabert-teissier(74):normal_cones_sheav_rel_jets} 
(in particular, from the finiteness proved in 
\cite[Theorem 4.15]{lejeune-jalabert-teissier(74):normal_cones_sheav_rel_jets}).

\begin{thm}\label{thm:multiplicity}
Let $\Y \subset \X$ be a closed subscheme, and let $m$ be a nonnegative integer. 
The set of $y \in \Y$ such that $\Y_{p(y)}$ has multiplicity at least $m$ at $y$ 
is Zariski-closed in $\Y$.
\end{thm}

For our purposes, $\Y\rightarrow B$ will be a family of subschemes of dimension 
$d$ in the family of (possibly reducible) varieties $\X\rightarrow B$. Since we are
interested in the multiplicities at several points at a time for the computation 
of Seshadri constants, we need a multipoint analogue of Theorem \ref{thm:multiplicity}.
Denote $\X_B^r:= \X \times_B \overset{r}{\cdots} \times_B \X \overset{p^r}{\rightarrow} B$ 
(respectively $\Y_B^r \rightarrow B$) the family whose fiber over $b$ is $(\X_b)^r$ 
(respectively $(\Y_b)^r$).  The next Corollary is an immediate consequence of \eqref{thm:multiplicity}.

\begin{cor}\label{cor:multiplicities}
  Let $\Y \subset \X$ be a closed subscheme, and let $m_1, \ldots, m_r$ be nonnegative integers. 
  The set of tuples 
  $(y_1, \ldots y_r) \in (\Y_b)^r$, $b \in B$,
  such that $\Y_b$ has multiplicity at least $m_i$ at each $y_i$ 
  is Zariski-closed in $\Y^r_B$.
\end{cor}

Given a closed subscheme $\Y \subset \X$ and a sequence $\m=(m_1, \ldots, m_r)$, denote
$\Sing_\m(\Y)\subset \Y^r_B$ the closed set given by \ref{cor:multiplicities}.
Fix be a relatively nef divisor $\LL$ on $\X$, and denote $\LL_b=\LL|_{\X_b}$.
Since the definition of multipoint Seshadri constants only makes sense
for distinct smooth points,
we will usually restrict to the complement of the diagonals and the
singularities of fibers in $\X^r_B$, which is Zariski open.

\begin{lem}\label{lem:semicontinuous}
  Let $r$ be a positive integer, let $\oover{\X}_B^r$ be the complement of the 
  diagonals and the singularities in $\X^r_B$, and let a real number $\epsilon >0$ 
  and a nonzero real vector ${\bf \ell}=(l_1,\cdots,l_r)$ be given with each $l_i\ge0$. 
  Then the set of tuples 
  $(p_1,\ldots,p_r)\in \oover{\X}_B^r$ such that 
  $$\epsilon_d(l_1p_1,\ldots,l_rp_r;\X_b,\LL_b)\le \epsilon$$ is the union of at most countably many 
  Zariski closed sets of $\oover{\X}_B^r$.
\end{lem}

\begin{proof}
  Let $\Hilb_d(\X/B)$ denote the relative Hilbert scheme of subschemes of $\X$ of
  (relative) dimension $d$.
  It has countably many irreducible components, which are irreducible projective schemes over $B$
  (see \cite{nitsure(05):constr_hilb_quot_sch}). Let $H$ be one of these components, and 
  let $\HH \subset \X\times_B H \rightarrow H$ be the corresponding universal family. 
  By the standard properties 
  of the Hilbert scheme, the intersection with $\LL_b$ is constant, i.e., $\LL_b^d\cdot \HH_{b,h}$ 
  does not depend on $b\in B$, $h\in H$. Denote this number by $\LL^d \cdot H$.

  For each choice of a sequence $\m=(m_1, \ldots, m_r)$ of $r$ nonnegative integers, 
  we apply corollary \ref{cor:multiplicities} above to the universal family 
  $\HH \rightarrow H$ and get a Zariski closed set 
  $$\Sing_\m(\HH)\subset\HH_H^r \subset (\X\times_B H)_H^r=\X_B^r \times_B H.$$
  Let $\USing_\m(\HH)\subset \X^r_B$ be the image of $\Sing_\m(\HH)$ by projection 
  on the first factor.
  Since $H$ is projective, $\USing_\m(\HH)$ is Zariski closed, and 
  $\oover{\USing}_\m(\HH):=\USing_\m(\HH)\cap\oover{\X}_B^r$ is Zariski closed in
  $\oover{\X}_B^r$. By the definition
  of Seshadri constants, for all $(p_1,\dots, p_r)\in \oover{\USing}_\m(\HH)$ 
  with $p_i \in \X_b$, $\epsilon_d(l_1p_1,\dots, l_rp_r;\X_b,\LL_b)\le \left(\LL^d \cdot H/\sum l_im_i\right)^{1/d}$.

  Since every $d$-dimensional subscheme is represented by a point in some component of the Hilbert scheme, it follows that
\begin{equation}
\label{eq:seshhilb} 
\epsilon_d(p_1,\dots, p_r;\X_b,\LL_b)=\inf_{(p_1,\dots,p_r) \in \oover{\USing}_\m(\HH)}
   \left\lbrace \left( \frac{L^d \cdot H}{\sum l_im_i} \right) ^{1/d}
  \right\rbrace.
\end{equation} 
  Thus, for each $\epsilon \in \R$, the set of tuples $(p_1,\ldots,p_r)$ such that 
  $\epsilon_d(p_1,\ldots,p_r;X,L)\le \epsilon$ is exactly
  $$\bigcup_{\{H: \frac{L^d \cdot H}{\sum l_im_i} \le \epsilon^d\}}  \oover{\USing}_\m(\HH),$$
  hence the claim.
\end{proof}

\begin{rmk}
\label{rmk:vgeneral} 
  From the previous Lemma \eqref{lem:semicontinuous} applied to $\X=X$, $B=\Spec K$,
  for very general points $p_1,\ldots,p_r$ the Seshadri constant
  $\epsilon_d(p_1,\ldots,p_r;X,L)$ and its weighted counterparts are independent of the 
  points chosen,
  and thus $\epsilon_d(r;X,L)$ is well defined. Moreover, (\ref{eq:seshhilb}) 
  shows that 
  $$\epsilon_d(r;X,L)=\inf_{\USing_\m(\HH)=X^r} 
   \left\lbrace \left( \frac{L^d \cdot H}{\sum m_i} \right)
     ^{1/d}\right\rbrace.$$
\end{rmk}
\begin{rmk}
  If $d=1$, then the infimum in (\ref{eq:seshhilb}) can be taken over
  all $\m$ and all components $\HH$ of the Hilbert scheme of $X$
  \emph{in all dimensions} (including, e.g., the isolated point
  corresponding to the whole variety $X$). Doing so,
  the Nakai-Moishezon for $\R$-divisors 
  \cite{campana-peternell(90):algeb_ample_cone_projec_variet} implies
  that the infimum is attained by some $\HH$ and $\m$ (see the proof
  of \cite[Proposition 4]{steffens(98):rem_sesh_consts}, and also
  \cite{debarre(04):seshad_const_abelian_variet},
  \cite{bauer-szemberg(01):local_pos_prin_pol_ab_3folds}) hence 
  the set of  values effectively taken by the Seshadri constant as
  points vary is either finite or countable.
\end{rmk}

\begin{rmk}
  For surfaces, a stronger version of Lemma \eqref{lem:semicontinuous} is known to hold.  Namely a finite number of Zariski closed sets (hence a single one) suffices,   and moreover the set of  values
  effectively taken by the Seshadri constant is either finite or has exactly one accumulation
  point which is $\sqrt{L^2/r}$. This has been proved by Oguiso \cite{oguiso(02):sesh_const_family_surf} 
  for $r=1$ and follows from Harbourne-Roé \cite{harbourne-roe(08):disc_beh_sesh_consts} for $r>1$. Unfortunately, the methods used to
  prove such finiteness do not seem to extend to varieties of
  higher dimension.
\end{rmk}

\begin{rmk}
  For Seshadri constants at higher dimensional centers (i.e. measuring
  the multiplicities at non-closed points) as defined by Paoletti in
  \cite{paoletti(94):sesh_consts_gonality_spc_curv_stable_lin_bdl}
  similar semicontinuity results hold; in this case, the Hilbert
  scheme has to be used as parameter space in the place of the $r$th
  product of the variety.
\end{rmk}
\section{Proofs of Theorems}

\subsection{Proof of Theorems \ref{thm:main} and \ref{thm:main2}}\label{sec:proofofmaintheorems}

We start with Theorem \ref{thm:main} and consider the degeneration to
the normal cone of very general points $p_1,\ldots, p_r$ in $X$.  In
detail set $B=\AA^1$ and let $\pi\colon \X\to X\times B$ be the blowup
at the points $\bar{p}_i=(p_i,0)\in X\times B$.  The exceptional
divisor $F$ is a disjoint union $F=\bigcup_i F_i$ where
$F_i=\PP(T_{p_i}\oplus \mathbb C)\simeq\PP^n$ is the projective
completion of the tangent space $T_{p_i}X$ of $X$ at $p_i$.  We denote
by $q_1,q_2$ the projections from $X\times B$ to the factors and let
$q=q_2\circ \pi\colon \X\to B$.

%
Set $\delta=\epsilon_d(r;X,L)$ and fix $1\le i\le r$.  By replacing $B$ with an open set around $0$ if
necessary we can find sections $\sigma'_{i,j}$ for $j=1,\ldots,s$
which pass through $p_i$ and have very general tangent direction
there.  Denote by $\sigma_{i,j}$ the section obtained from the proper
transform of $\sigma'_{i,j}$ in $\X$.  Then the collection
$\S_i=\{\sigma_{i,1}(0),\ldots,\sigma_{i,s}(0)\}$ is a set of $s$
very general points in $F_i\simeq \PP^n$.  
Denote $\sigma':=(\sigma'_{1,1}, \dots, \sigma'_{r,s}):B \rightarrow X^{rs}$
and $\sigma:=(\sigma_{1,1}, \dots, \sigma_{r,s}):B \rightarrow \X_B^{rs}$.

Each component of the relative Hilbert scheme of $X\times B$ is of the
form $H'=Z \times B$ where $Z$ is a component of the Hilbert scheme of
$X$; for each such $H'$ there is a unique 
component $H$ of the Hilbert scheme of $\X$ with 
$H' \cap q_2^{-1}(B\setminus\{0\})=H \cap q^{-1}(B\setminus\{0\})$. 
We resume notations as in the previous section, with 
$\Z \rightarrow Z$, $\HH' \rightarrow H'$, $\HH \rightarrow H$ as universal Hilbert families on 
$X$, $X \times B$ and $\X$ respectively. 

Recyling notation from above, let $\USing_\m(Z)$ be the projection of
$\Sing_{\m}(Z)$ to $X^{rs}$ and consider for a while a sequence
$\m=(m_1, \ldots, m_{rs})$ of $rs$ 
nonnegative integers such that $\USing_\m(\Z)=X^{rs}$. For
convenience, set also $m_{i,j}=m_{(i-1)s+j}$ for $i=1, \ldots, r$,
  $j=1, \ldots, s$. Then
$\USing_\m(\HH')=(X\times B)_B^{rs}=X^{rs} \times B$ and $\USing_\m(\HH)=\X_B^{rs}$.
Consider the following pullback diagram:
$$
\xymatrix{
\Sing_\m(\HH')\times_{\USing_\m(\HH')} B \ar[d]^{\tau_1} \ar[r]^{\tau_2} 
& B\ar[d]^{\sigma'\times i}\\
 \Sing_\m(\HH') \ar@{->>}[r] & \USing_\m(\HH')=X^{rs} \times B}
$$
$\tau_2$ is onto and proper, so (restricting to a smaller neighbourhood of $0$ 
if needed) we may choose a section; denote by $\zeta:B \rightarrow Z$ the 
composition of this section with the natural maps 
$$\Sing_\m(\HH')\times_{\USing_\m(\HH')} B
\overset{\tau_1}{\rightarrow} \Sing_\m(\HH')\hookrightarrow {\HH'}_{H'}^{rs} \hookrightarrow X^{rs}\times H
\rightarrow H=Z \times B \rightarrow Z.$$
Let $\Y'\rightarrow B$ be the family obtained from the universal family $\Z \rightarrow Z$ 
by base change through $\zeta$. By construction, each fiber $\Y_b'$ with $b\ne 0$ 
is a subscheme of $X$ of dimension $d$ with 
a point of multiplicity $\ge m_{i,j}$ at $\sigma_{i,j}(b)$.
Consider the strict transform $\Y$ of $\Y'$ in $\X$. By flatness, $\Y_0$ has dimension $d$, and by
semicontinuity of multiplicities (\ref{cor:multiplicities}) it has 
a point of multiplicity $\ge m_{i,j}$ at $\sigma_{i,j}(0)$. Therefore for every
$i$ such that some $m_{i,j}>0$, $\Y_0 \cap F_i$ is a $d$-dimensional subscheme of 
$F_i \cong \PP^n$ of
degree at least $m_i:=(\epsilon_d(s;\PP^n,\OO_{\PP^n}(1)))^d \sum_j m_{i,j}$. Since this degree
is exactly the multiplicity of $\Y'_0$ at $p_i$, it follows that
$$L \cdot Z = L \cdot \Y'_0 \ge \delta^d\sum_{i=1}^r m_i =
\delta^d (\epsilon_d(s;\PP^n,\OO_{\PP^n}(1)))^d \sum_{k=1}^{rs} m_{k}.$$
Since this is true whenever $\USing_\m(\Z)=X^{rs}$, in view of \ref{rmk:vgeneral} 
the claimed bound on $\epsilon_d(rs;X,L)$ follows. 

The proof of Theorem \ref{thm:main2} follows easily as if $s_1\le s_2\le \cdots \le s_r$ then $s=\sum_{i=1}^r s_i \le rs_r$ and
$$\epsilon_d(s;X,L) \ge \epsilon_d(rs_r;X,L)\ge \epsilon_d(r;X,L)\cdot
\epsilon_d(s_r;\PP^n,\OO_{\PP^n}(1)).\qed $$

\subsection{Proof of Theorem
  \ref{thm:glueing}}\label{sec:proofofthemglueing}
Essentially, the proof of Theorem \ref{thm:main} works also in this
more general setting. We have proved the particular case first for clarity,
and give next a sketch of the changes needed for \ref{thm:glueing}.

Just as above let $\X\to B$ be the degeneration to the normal cone of
very general points $p_1,\ldots, p_r$ in $X$, with exceptional divisors $\PP^n$.
Also just as above, by shrinking $B$ if necessary, for every component $Z$ of 
$\Hilb_d(X)$ and every $\m=(m_1,\dots,m_{s})$ with $\USing_\m(\Z)=X^{s}$
there exist schemes in $\Hilb_d(X)$ with multiplicities at least
$$m_i:=(\epsilon_d(s_i,{\bf \ell}^i;\PP^n,\OO_{\PP^n}(1)))^d 
\sum_{j=s_1+\cdots+s_{i-1}+1}^{s_1+\cdots+s_{i}} m_j$$
at general points $p_i$, $i=1, \ldots, r$. From this the result
follows exactly as in the previous case. \qed

\bibliography{biblio2}

{\small \noindent {\tt j.ross@dpmms.cam.ac.uk}} \newline
\noindent DPMMS, Centre for Mathematical Sciences, Wilberforce Road, Cambridge CB3 0WA. UNITED KINGDOM. 

{\small \noindent {\tt jroe@mat.uab.cat}} \newline
\noindent Departament de Matemàtiques, Universitat Aubònoma de
Barcelona. 08193 Bellaterra (Barcelona). SPAIN. \\

\end{document}